\documentclass[12pt,a4paper,leqno]{article}

\usepackage{latexsym}
\usepackage{amssymb,exscale}
\usepackage[centertags]{amsmath}
\usepackage{amsthm}
\usepackage{graphicx}

\usepackage[dvips]{hyperref}

\theoremstyle{definition}
\newtheorem{lemma}{Lemma}
\newtheorem{theorem}{Theorem}
\newtheorem{remark}{Remark}
\newtheorem{question}{Question}

\renewcommand{\phi}{\varphi}

\renewcommand{\(}{\bigl(}
\renewcommand{\)}{\bigr)\vphantom{)}}

\renewcommand{\div}{\operatorname{div}}

\newcommand{\Ex}{\mathbb E\,}
\newcommand{\R}{\mathbb R}

\newcommand{\Z}{\mathbb Z}

\newcommand{\almost}[1]{$#1$\nobreakdash-\hspace{0pt}almost}

\renewcommand{\Pr}[1]{\mathbb{P}\mskip1.5mu\(\mskip1.5mu#1\mskip1.5mu\)}

\begin{document}

\title{Divergence of a stationary random vector field can be always
  positive (a Weiss' phenomenon)}

\author{Boris Tsirelson}

\date{}
\maketitle

\stepcounter{footnote}
\footnotetext{%
 This research was supported by \textsc{the israel science foundation}
 (grant No.~683/05).}

\begin{abstract}
The divergence of a stationary random vector field at a given point is
usually a centered (that is, zero mean) random variable. Strangely
enough, it can be equal to $ 1 $ almost surely. This fact is another
form of a phenomenon disclosed by B.~Weiss in 1997.
\end{abstract}

\section*{Introduction}
If a random vector field is stationary (that is, shift-invariant in
distribution), then the expectation of its divergence must vanish,
since it is equal to the divergence of the expectation, thus, the
divergence of a constant vector field. This simple argument is
conclusive provided that the expectations are well-defined
(especially, for Gaussian processes). Waiving existence of first
moments we cannot ask about the expected divergence, but we still may
ask, whether the divergence can be always (strictly) positive, or
not. The answer appears to be negative in dimension one but
affirmative in dimension two. The former is evident, while the latter
is demonstrated by a non-evident counterexample constructed below
following an idea of B.~Weiss \cite{Weiss}.

The phenomena under consideration manifest themselves equally well in
two forms, discrete and continuous. Starting with the discrete setup
we treat two lattices $ \Z $ and $ \Z^2 $ as two (infinite) graphs.
\[
\begin{gathered}\includegraphics{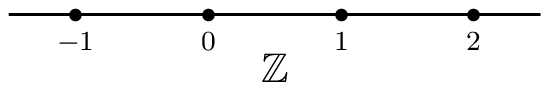}\end{gathered}\qquad
\begin{gathered}\includegraphics{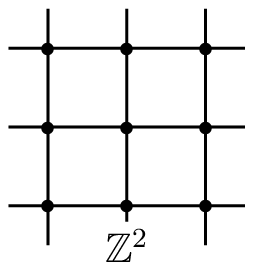}\end{gathered}
\]
By a vector field me mean a real-valued function on oriented edges
such that its sum over the two orientations of an edge vanishes. By
the divergence of a vector field we mean the following real-valued
function on vertices: given a vertex, we sum over the (two or four)
outgoing edges the values of the vector field.
\[
\includegraphics{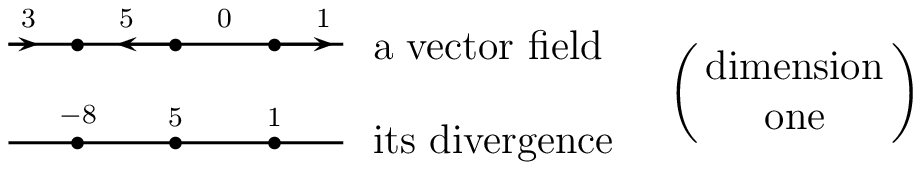}
\]
For a vector field $ v = (v_{x,x+1})_{x\in\Z} $ on $ \Z $, its
divergence is given by
\[
(\div v)_x = v_{x,x+1} - v_{x-1,x} \quad \text{for } x \in \Z \, .
\]
If a random vector field $ v $ on $ \Z $ is stationary then the two
random variables $ v_{-1,0} $ and $ v_{0,1} $ are identically
distributed. If also $ (\div v)_0 \ge 0 $ a.s., then $ v_{-1,0} \le
v_{0,1} $ a.s., therefore $ \Pr{ v_{-1,0} \le a < v_{0,1} } = \Pr{
v_{-1,0} \le a } - \Pr{ v_{0,1} \le a } = 0 $ for all $ a \in \R
$. Letting $ a $ run over a dense countable set we get $ v_{-1,0} =
v_{0,1} $ a.s. It means that
\begin{quotation}
\noindent a stationary random vector field $ v $ on $ \Z $ cannot
satisfy the condition $ \Pr{ (\div v)_0 > 0 } = 1 $.
\end{quotation}
In dimension two the situation is different.

\begin{theorem}
There exists a stationary random vector field $ v $ on $ \Z^2 $ whose
divergence is equal to $ 1 $ everywhere, almost sure.
\end{theorem}

\noindent (For the proof see Sect.~\ref{sect1}.) The corresponding
result in the continuous setup follows by smoothing, namely,
convolution with the indicator function of the square $ (-0.5,0.5)
\times (-0.5,0.5) $. In other words: if, say, $ v=1 $ on the edge $
((0,0),(0,1)) $ and $ v=0 $ on other edges of $ \Z^2 $, then the first
(horizontal) component of the smoothed vector field on $ \R^2 $ at $
(x,y) $ is equal to $ 1-|x-0.5| $ if $ |x-0.5|<1 $, $ |y|<0.5 $ (and $
0 $ otherwise).
\begin{gather*}
\text{vector field}\\
\begin{gathered}\includegraphics{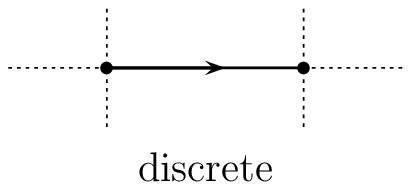}\end{gathered}\qquad
\begin{gathered}\includegraphics{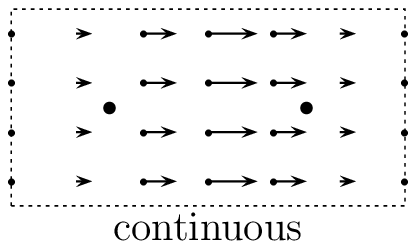}\end{gathered}
\end{gather*}
The divergence of the smoothed field is equal to $ +1 $ on the square $
(-0.5,0.5) \times (-0.5,0.5) $ and $ -1 $ on $ (0.5,1.5) \times
(-0.5,0.5) $; just the smoothed $ \div v $.
\begin{gather*}
\text{divergence}\\
\begin{gathered}\includegraphics{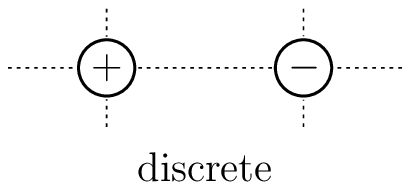}\end{gathered}\qquad
\begin{gathered}\includegraphics{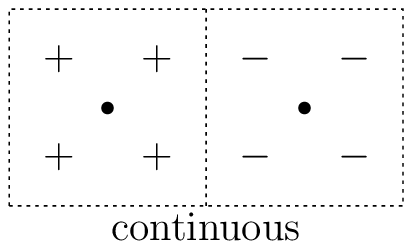}\end{gathered}
\end{gather*}
Having $ v $ such that $ \div v = 1 $ everywhere on $ \Z^2 $ we get
the smoothed divergence equal to $ 1 $ almost everywhere on $ \R^2 $
(and in the distributional sense). Some additional smoothing gives a
smooth vector field of divergence $ 1 $ everywhere.

\section[The construction and the proof]{The construction and the
 proof}
\label{sect1}
In order to keep the matter as discrete as possible, from now on we
consider only \emph{integer-valued} vector fields on $ \Z^2 $.

\begin{lemma}\label{lemma1}
If there exist stationary random vector fields $ v_1,v_2,\dots $ such
that

(a) $ \Pr{ (\div v_n)_x = 1 } \to 1 $ as $ n \to \infty $, for every
vertex $ x $ of the graph $ \Z^2 $,

(b) $ \sup_n \Pr{ | (v_n)_y | > C } \to 0 $ as $ C \to \infty $, for
every edge $ y $ of the graph~$ \Z^2 $,

\noindent then there exists a stationary random vector field $ v $
such that

(A) $ \Pr{ (\div v)_x = 1 } = 1 $ for every vertex $ x $ of the graph
$ \Z^2 $,

(B) $ \Pr{ | v_y | > C } \le \sup_n \Pr{ | (v_n)_y | > C } $ for
every $ C $ and every edge $ y $ of the graph $ \Z^2 $.
\end{lemma}

\begin{proof}
The distribution $ \mu_n $ of $ v_n $ is a probability measure on the
space $ \Z^E $ of all maps $ E \to \Z $; here $ E $ is the (countable)
set of all edges of $ \Z^2 $. Using the one-point compactification $
\overline\Z = \Z \cup \{\infty\} $ of $ \Z $ we may treat $ \mu_n $ as
measures on the compact metrizable space $ \overline\Z^E $. All
probability measures on $ \overline\Z^E $ being a compact metrizable
space, we take a convergent subsequence: $ \mu_{n_k} \to \mu $.

Let $ v $ be distributed $ \mu $, then (B) is satisfied, and all
values of $ v $ are finite a.s.\ due to (b) and (B). Clearly, $ v $ is
stationary. Treating $ (\div v)_x $ as a function (of $ v $) defined
and continuous \almost{\mu} everywhere (namely, on $ \Z^E $) we see
that $ (\div v_{n_k})_x $ converges in distribution to $ (\div v)_x $
as $ k \to \infty $. Thus, (a) implies (A).
\end{proof}

Random vector fields $ v_n $ will be constructed out of non-random
finite fragments. For example, $ n=2 $:
\[
\begin{gathered}\includegraphics{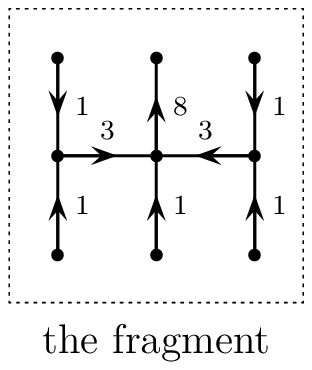}\end{gathered}\qquad
\begin{gathered}\includegraphics{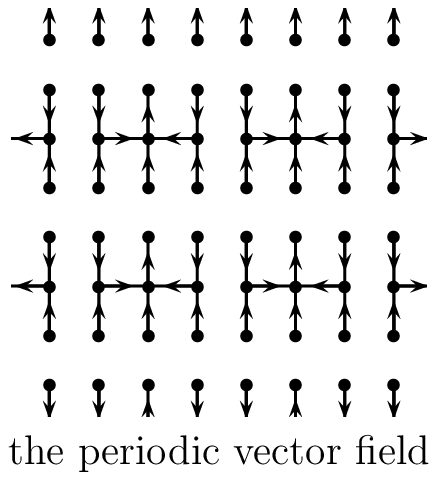}\end{gathered}
\]
The fragment of size $ 3 \times 3 $ (whose construction will be
explained later) is repeated, forming a (double-) periodic vector
field (non-random). Shifting this periodic field we get $ 3 \cdot 3 =
9 $ periodic fields. These $ 9 $ vector fields are the possible values
of the random vector field $ v_2 $; they have equal probabilities ($
1/9 $), by definition (of $ v_2 $). Thus, $ v_2 $ is stationary.

Note that the divergence of $ v_2 $ at a given point (say, the origin)
takes on two values, $ +1 $ (with probability $ 8/9 $) and $ -8 $
(with probability $ 1/9 $). Also, the value of $ v_2 $ on a given
horizontal edge takes on three values $ -3, 0, 3 $ with probabilities
$ 1/9, 7/9, 1/9 $ respectively.

The fragment of size $ 3 \times 3 $, used above, is the second element
of a sequence, whose $ n $-th element is of size $ (2^n-1) \times
(2^n-1) $. The sequence is constructed recursively. Its first element,
of size $ 1 \times 1 $, is trivial: just a single vertex, no
edges. The $ (n+1) $-th fragment contains four copies of the $ n $-th
fragment (two of them being turned upside down) connected as follows:
\[
\begin{gathered}\includegraphics{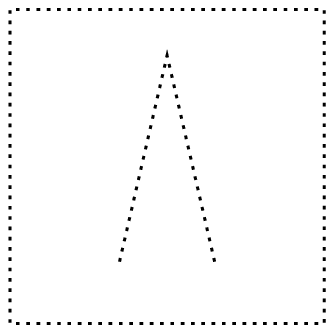}\end{gathered}
\quad = \quad
\begin{gathered}\includegraphics{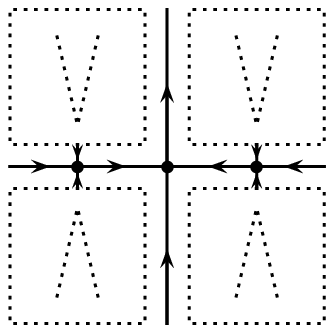}\end{gathered}
\]
Here are the first three fragments:
\[
\begin{gathered}\includegraphics{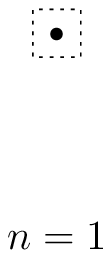}\end{gathered}
\begin{gathered}\includegraphics{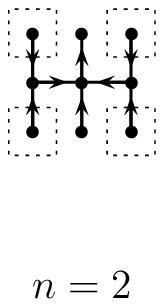}\end{gathered}
\quad
\begin{gathered}\includegraphics{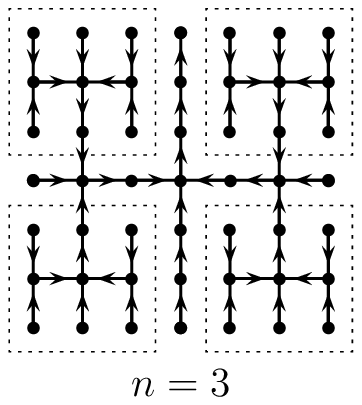}\end{gathered}
\]
For each fragment, the vector field conforms to the given oriented
graph, and its divergence is equal to $ +1 $ at all vertices except
for one vertex. Such vector field exists and is unique, since the
graph is a tree; the divergence at the root is equal to $ - \(
(2^n-1)^2 - 1 \) = - ( 2^{2n} - 2^{n+1} ) $. Note the consistency: the
four copies of the $ n $-th fragment occur in the $ (n+1) $-th
fragment as vector fields (not only graphs).

For each $ n $ we construct $ v_n $ out of the $ n $-th fragment in
the same way as we did it for $ n=2 $. Clearly, $ v_n $ is an
integer-valued stationary random vector field on $ \Z^2 $, and
\[
\Pr{ (\div v_n)_x = 1 } = 1 - \frac1{ (2^n-1)^2 } \, ,
\]
which verifies Condition (a) of Lemma \ref{lemma1}. It remains to
check Condition (b), which is the point of the next lemma.

\begin{lemma}\label{lemma2}
$ \sup_n \Pr{ | (v_n)_y | > C } = O(1/\sqrt C) $ as $ C \to \infty $,
for every edge $ y $ of the graph $ \Z^2 $.
\end{lemma}

\begin{proof}
First, the maximal possible value of $ | (v_n)_y | $ is equal to $
2^{2n} - 2^{n+1} $ (this is the value on the edge that enters the root
of the tree). Second,
\[
\Pr{ | (v_{n+1})_y | \le 2^{2n} - 2^{n+1} } \ge \frac{ 4 \cdot
(2^n-1)^2 }{ (2^{n+1}-1)^2 }
\]
since the $ (n+1) $-th fragment has $ (2^{n+1}-1)^2 $ vertices, and
out of these, $ 4 \cdot (2^n-1)^2 $ vertices belong to the four copies
of the $ n $-th fragment. (For edges the ratio is even larger.)
Similarly,
\[
\Pr{ | (v_{n+i})_y | \le 2^{2n} - 2^{n+1} } \ge \frac{ 4^i \cdot
(2^n-1)^2 }{ (2^{n+i}-1)^2 }
\]
for $ i = 1,2,\dots\, $ Therefore
\begin{multline*}
\Pr{ | (v_{n+i})_y | \le 2^{2n} } \ge \bigg( \frac{ 2^{n+i} - 2^i }{
 2^{n+i} - 1 } \bigg)^2 = \bigg( 1 - \frac{ 2^i-1 }{ 2^{n+i} - 1 }
 \bigg)^2 \ge (1-2^{-n})^2 \ge \\
 \ge 1 - 2 \cdot 2^{-n} \, ;
\end{multline*}
\[
\sup_n \Pr{ | (v_n)_y | > 2^{2k} } \le 2 \cdot 2^{-k} \, ;
\]
the lemma follows.
\end{proof}

\section[Remarks and questions]{Remarks and questions}
\label{sect2}
\begin{remark}
It is easy to see that the constructed sequence $ (v_n)_n $ converges
in distribution. No need to choose a subsequence (as in the proof of
Lemma~\ref{lemma1}).
\end{remark}

\begin{question}
Is the condition $ \Ex \sqrt{ |v_y| } < \infty $ compatible with $
\div v = 1 $?
\end{question}

\begin{remark}
As shown by B.~Weiss \cite{Weiss}, the sample functions of a
stationary complex-valued process on the complex plane can be
non-constant entire functions.
\end{remark}

\bigskip
\filbreak
{
\small
\begin{sc}
\parindent=0pt\baselineskip=12pt
\parbox{4in}{
Boris Tsirelson\\
School of Mathematics\\
Tel Aviv University\\
Tel Aviv 69978, Israel
\smallskip
\par\quad\href{mailto:tsirel@post.tau.ac.il}{\tt
 mailto:tsirel@post.tau.ac.il}
\par\quad\href{http://www.tau.ac.il/~tsirel/}{\tt
 http://www.tau.ac.il/\textasciitilde tsirel/}
}

\end{sc}
}
\filbreak

\end{document}